\documentclass{amsart}
\usepackage{amssymb} 
\usepackage{color} 
\newtheorem{corollary}{Corollary}
\newtheorem{theorem}{Theorem}
\newtheorem{lemma}{Lemma}
\newtheorem{proposition}{Proposition}

\begin{document}
\title{Pauli gradings on Lie superalgebras and graded codimension growth}

\author[D. D. Repov\v s and M. V. Zaicev]
{Du\v san D. Repov\v s and Mikhail V. Zaicev}

\address{Du\v san D. Repov\v s \\Faculty of Education, and
Faculty of  Mathematics and Physics, University of Ljubljana,
 Ljubljana, 1000, Slovenia}
\email{dusan.repovs@guest.arnes.si}

\address{Mikhail V. Zaicev \\Department of Algebra\\ Faculty of Mathematics and
Mechanics\\  Moscow State University \\ Moscow,119992, Russia}

\email{zaicevmv@mail.ru}

\thanks{The first author was supported by the Slovenian Research Agency
grants P1-0292-0101, J1-7025-0101, J1-6721-0101 and J1-5435-0101. 
The second author was supported 
by the Russian Science Foundation, grant 16-11-10013.
We thank the reviewers for comments.}



\begin{abstract}
We introduce grading on certain finite dimensional simple Lie superalgebras
of type $P(t)$ by elementary abelian 2-group. This grading gives rise to Pauli matrices
and is a far generalization of $(\mathbb Z_2\times \mathbb Z_2)$-grading on Lie algebra
of $(2\times 2)$-traceless matrices.We use this grading for studying numerical invariants 
of polyomial identities of Lie superalgebras. In particular, we compute graded PI-exponent
corresponding to Pauli grading.
\end{abstract}

\date{\today}

\maketitle
\vskip 0.2in
%{\Large

\section{introduction}

In this paper we study algebras over a field $F$ of characteristic zero. Group graded algebras
have been intensively studied in the last decades (see, for example, 
\cite{BSZ2001, BShZ, BZ,  Eld, EldKo, NVanO, PZ1988, SegZai}). 
All possible gradings on
matrix algebras over an algebraically closed field were described in \cite{BSZ2001, BZ}.
Recently, all gradings by a finite abelian groups on finite dimensional simple real algebras
have also been
classified in \cite{BZ1, R-E}. 
Many authors have also paid attention to  grading 
on Lie algebras 
\cite{BShZ, DrEld, EldKo, PZ1988}. 
Both, in associative and Lie 
case, an exceptional role is played by gradings which cannot be "refined" --
in particular, gradings
whose homogeneous components are one-dimensional 
\cite{ BSZ2001, BZ, DrEld, PZ1988}.
Classification of group gradings on Lie superalgebras is only in its initial stages (see, e.g., \cite{B-ArX}). 
Therefore an important role is played by new examples of gradings on Lie superalgebras.

It is well known that abelian gradings are closely connected to
 automorphism and involution actions 
on algebra (see, for example, \cite{BSZ2001}), 
hence the knowledge of gradings gives us an important information about the group 
of automorphisms and antiautomorphisms of an algebra. 
Another application of gradings is the study 
of graded and non-graded identities and their numerical invariants.

Given an algebra $A$, one can associate to it an infinite sequence of non-negative integers 
$$\{c_n(A)\}, \ \ 
n=1,2,\ldots,\,$$
 called {\it codimensions} of $A$. 
The study of asymptotic behavior of $\{c_n(A)\}$ is
one of the most important and current approaches in the modern PI-theory \cite{GZbook}. 
In many cases
codimension growth is exponentially bounded. 
In particular, 
$$\dim A=d<\infty \ \ \Rightarrow \ \   c_n(A)\le d^{n+1}$$
(see \cite{BDLaa} and also \cite[Proposition 2]{GZTams2010}). 
If, in addition, $A$ is endowed with a 
grading by a group $G$ then one can
also define the graded codimension sequence $c_n^G(A)$. 
For a finite dimensional algebra $A,$ graded 
and ordinary codimensions satisfy the following inequalities:
\begin{equation}\label{eq1}
c_n(A)\le c_n^G(A)\le (\dim A)^{n+1}
\end{equation}
(see  \cite{BDLaa}).

As a rule, an investigation of asymptotics of graded codimensions is much easier 
than a study of non-graded codimensions. This fact was used in our previous papers for obtaining
the results on  both graded and non-graded codimension growth 
\cite{GZ-JLMS, RZ1, RZ2, RZ3}.

If $A$ is a finite dimensional graded simple algebra then there exist the limits
\begin{equation}\label{eq1a}
exp(A)=\lim_{n\to\infty} \sqrt[n]{c_n(A)},\quad
exp^{G}(A)=\lim_{n\to\infty} \sqrt[n]{c_n^{G}(A)}
\end{equation}
and according to (\ref{eq1}) we have
\begin{equation}\label{eq2}
exp(A)\le exp^G(A)\le \dim A.
\end{equation}

It is well known that in many most important cases of algebras (associative, Lie, Jordan, alternative, etc.)
\begin{equation}\label{eq3}
exp(A)= \dim A,
\end{equation}
provided that $A$ is simple and $F$ is  algebraically closed
 \cite{GShZ, GZ1998, Z}.
In this case $exp^G(A)$ is also equal to $\dim A$ for any grading on $A$. If $A$ is graded simple but
not simple in the usual sense then graded and non-graded exponents can
 differ. For example, if $G$ is a
finite abelian group of order $|G|=m$ and $A$ is its  group algebra, $A=FG$, then $exp(A)=1$
whereas $exp^G(A)=m$. Clearly, if $A$ is simple in non-graded sense then $A$ is also graded simple for 
any $G$-grading. Relations (\ref{eq2}) and
 (\ref{eq3}) show that the conjecture that $exp(A)=exp^G(A)$ holds for associative, 
Lie, Jordan and alternative algebras over an algebraically closed field.

Nevertheless, in the Lie superalgebra case there exist simple algebras such that $exp(A)$
and
 $exp^G(A)$ exist and 
are strictly less than $\dim A$ (see \cite{GZ-JLMS, RZ3}).
 Here we are talking  about canonical $\mathbb Z_2$-grading
on Lie superalgebras. Therefore the study of relations between graded and non-graded PI-exponents
is of interest in the general case. In particular, if the conjecture that $exp(A)=exp^G(A)$ is confirmed
then it would give us a powerful tool for computing precise asymptotics of codimension growth. Another
consequence would be the independence of $exp^G(A)$ on the particular $G$-grading.

The goal of the present paper is twofold. In the first part  we define the so-called Pauli $G$-grading
on the simple Lie superalgebra of the type $L=P(t)$ (in the notation of \cite{Kac}, for general material on Lie superalgebras see alo \cite{Sch}),
 where $t$ is 
the power of 2 and $G$ is an elementary abelian 2-group. 
This grading posesses many remarkable properties. 
In fact, it is induced from the grading on simple 3-dimensional Lie algebra $sl_2(F)$ by Pauli matrices and is compatible with the canonical $\mathbb Z_2$-grading. All non-zero homogeneous components of $L$ are 
one-dimensional. Also, any even homogeneous element $0\ne a\in L_g$ is a non-degenerate matrix and for any
homogeneous elements $a\in L_g, b\in L_h$ their Lie supercommutator is either zero or non-degenerate.
In the second part of the paper we investigate the graded codimension growth of $L$. We show that all 
computations are much easier than in the non-graded case due to the remarkable properties of Pauli grading.

Our main result is Theorem \ref{T1} below, stating that $exp^G(P(t))=t^2-1+t\sqrt{t^2-1}$.
Note that
Theorem \ref{T1} 
is true for $t=2$ although $P(2)$ is not simple and $exp^G(P(2))=3+2\sqrt 3$ holds for both
Pauli grading and the canonical $\mathbb Z_2$-grading (see \cite{RZ1}).

\begin{theorem}\label{T1}
Let $L$ be a Lie superalgebra of the type $P(t)$, $t=2^q, q\ge 1$, equipped with $G$-grading
given in Proposition \ref{P2}. Then $G$-graded PI-exponent of $L$ exists and
$$
exp^G(L)=t^2-1+t\sqrt{t^2-1}.
$$
\end{theorem}

\section{Pauli gradings}

Let $L$ be an algebra over a field $F$ and let $G$ be a group. One says that $L$ is $G$-graded if $L$
has a vector space decomposition
$$
L=\bigoplus_{g\in G}L_g
$$
such that $L_gL_h\subseteq L_{gh}$ for all $g,h\in G$. Subspaces $L_g, g\in G$, are called homogeneous 
components of  $L$. Any element $a\in L_g$ is called homogeneous of degree $\deg a=g$. The subset
$$
Supp~L=\{g\in G|L_g\ne 0\}
$$
is said to be the support of the grading. A subspace $V\subseteq L$ is called homogeneous if
$$V=\bigoplus_{g\in G}V\cap L_g.$$

Let $A$ and $B$ be two associative algebras and let $G$ and $H$ be two groups. Suppose that $A$ and 
$B$ are endowed by $G$- and $H$-gradings, respectively,
$$
A=\bigoplus_{g\in G}A_g, \quad B=\bigoplus_{h\in H}B_h.
$$
Then one can introduce $G\times H$-grading on the tensor product $A\otimes B$ by setting
$$
(A\otimes B)_{gh}=A_g\otimes B_h.
$$

An associative algebra $R$ is said to be a superalgebra if $R$ has  some $\mathbb Z_2$-grading, that is
$$
R=R^{(0)}\oplus R^{(1)},~R^{(0)}R^{(0)}+R^{(1)}R^{(1)}\subseteq R^{(0)},~
R^{(0)}R^{(1)}+R^{(1)}R^{(0}\subseteq R^{(1)}.
$$
A special case of associative superalgebras which we will use later is the $\mathbb Z_2$-graded
$n\times n$ matrix algebra $R=M_{k,l}(F)$ with
$$
R= \left\{ \left(
           \begin{array}{cc}
             A & B \\
             C & D \\
           \end{array}
         \right)  \right\}=R^{(0)}\oplus R^{(1)}
         \, ,\quad
         R^{(0)} = \left\{ \left(
           \begin{array}{cc}
             A & 0 \\
             0 & D \\
           \end{array}
         \right)\right\},\,
               R^{(1)} = \left\{ \left(
           \begin{array}{cc}
             0 & B \\
             C & 0 \\
           \end{array}
         \right)
\right\}
$$      
where $n=k+l$, $A,B,C,D$ are $k\times k,k\times l, l\times k$ and $l\times l$  matrices,
respectively. In particular, when $k=l$ we have $\mathbb Z_2$-grading on $M_{2k}(F)$
which will be used for the definition of Lie superalgebra $P(k)$.

Recall now that $\mathbb Z_2$-graded non-associative algebra $L=L^{(0)}\oplus L^{(1)}$ is
called a {\it Lie superalgebra} if it satisfies homogeneous relations
$$
ab+(-1)^{|a||b|}ba=0,~a(bc)=(ab)c+(-1)^{|a||b|}b(ac)=0
$$
for all $a,b,c\in L^{(0)}\cup L^{(1)}$ where $|x|=0$ if $x\in L^{(0)}$ and $|x|=1$ if $x\in L^{(1)}$.
In particular, any associative superalgebra $R=R^{(0)}\oplus R^{(1)}$ with the new product called
supercommutator, defined for homogeneous elements  as
$$
[a,b]= ab-(-1)^{|a||b|}ba
$$
becomes a Lie superalgebra.

Let $L^{(0)}\oplus L^{(1)}$ be a Lie superalgebra and let $G$ be a group. Then a $G$-grading 
$$L=\oplus_{g\in G}L_g$$ is called {\it compatible} with $\mathbb Z_2$-grading of $L$ if
$L_g\subseteq L^{(0)}$ or $L_g\subseteq L^{(1)}$ for all $g\in G$.

For defining the Pauli grading on the associative matrix algebra $M_{2^q}(F)$ we start with $q=1$.
Consider  $2\times 2$ matrices
\begin{equation}\label{e*}
\sigma_0= \left\{ \left(
           \begin{array}{cc}
             1 & 0 \\
             0 & 1 \\
           \end{array}
         \right)  \right\},
\sigma_1= \left\{ \left(
           \begin{array}{cc}
             1 & 0 \\
             0 & -1 \\
           \end{array}
         \right)  \right\},
\sigma_2= \left\{ \left(
           \begin{array}{cc}
             0 & 1 \\
             1 & 0 \\
           \end{array}
         \right)  \right\},
\sigma_3= \left\{ \left(
           \begin{array}{cc}
             0 & 1 \\
             -1 & 0 \\
           \end{array}
         \right)  \right\}.
\end{equation}
Matrices (\ref{e*}) are closely related to Pauli matrices.
$$
\sigma_x= \left\{ \left(
           \begin{array}{cc}
             0 & 1 \\
             1 & 0 \\
           \end{array}
         \right)  \right\},
\sigma_y= \left\{ \left(
           \begin{array}{cc}
             0 & -i \\
             i & 0 \\
           \end{array}
         \right)  \right\},
\sigma_z= \left\{ \left(
           \begin{array}{cc}
             1 & 0 \\
             0 & -1 \\
           \end{array}
         \right)  \right\}.
$$
It is well-known that the linear span $L=<\sigma_x,\sigma_y,\sigma_z>$ is closed under Lie
commutator and $L\simeq su(2)$ as Lie algebra whereas the span $<\sigma_0,\sigma_1,\sigma_2,
\sigma_3>$ as an associative algebra is isomorphic to $M_2(F)$. Denote by $G=<a>_2\times <b>_2$
the product of two cyclic groups of order 2 with generators $a$ and $b$, respectively. Clearly,
$G$ is isomorphic to $\mathbb Z_2\times \mathbb Z_2$ and the decomposition
\begin{equation}\label{g1}
R=M_2(F)=R_e\oplus R_a\oplus R_b\oplus R_{ab}
\end{equation}
is a $G$-grading, where
$$
R_e=<\sigma_0>, R_a=<\sigma_1>,R_b=<\sigma_2>,R_{ab}=<\sigma_3>.
$$
We call the grading (\ref{g1}) on $M_2(F)$ {\it Pauli grading} on $M_{2}(F)$.

We generalize this construction to matrices of arbirary size $2^q, q\ge 2$ in the following way. Let
$R=R_1\otimes\cdots\otimes R_q$ where all $R_1,\ldots, R_q$ are isomorphic to the $2\times 2$ matrix
algebra $M_2(F)$. Let also 
\begin{equation}\label{al0}
G_0=G_1\times\ldots\times G_q, G_j=<a_j>_2\times<b_j>_2\simeq \mathbb Z_2\oplus \mathbb Z_2,j=1,\ldots,q.
\end{equation}
Then $R$ has a basis consisting of elements
\begin{equation}\label{e**}
c=x_1\otimes\cdots\otimes x_q
\end{equation}
where all $x_1,\ldots,x_q$ are of the type (\ref{e*}). Then in the Kronecker realization of
tensor product of matrices for transpose involution $T$ we have
$$
c^T=(x_1\otimes\cdots\otimes x_q)^T=x_1^T\otimes\cdots\otimes x_q^T.
$$
In particular, the element $c$ of the type (\ref{e*}) is symmetric if and only if the number
of matrices $\sigma_3$ among $x_1,\ldots,x_q$ is even and $c^T=-c$ if and only if the number
of $\sigma_3$ is odd.

All $R_1,\ldots, R_q$ have Pauli grading as defined earlier and we can extend these gradings to 
their tensor product $R$. Then we obtain $G_0$-grading on $R$
$$
R=\bigoplus_{g\in G_0} R_g
$$
where $R_g=<x_1\otimes\cdots\otimes x_q>$ and all $x_1,\ldots,x_q$ are of the type (\ref{e*}).
Moreover, we have
\begin{equation}\label{al1}
\deg(x_1\otimes\cdots\otimes x_q)=\deg x_1\cdots \deg x_q
\end{equation}
where
\begin{equation}\label{al2}
\deg x_i=
\left\{
  \begin{array}{rcl}
     e_i, &\quad \hbox{if} \quad & x_i=\sigma_0 \\
    a_i, &\quad \hbox{if} \quad & x_i=\sigma_1\, ,\\
    b_i, &\quad \hbox{if} \quad & x_i=\sigma_2\, ,\\
    a_ib_i, &\quad \hbox{if} \quad & x_i=\sigma_3\, \\
               \end{array}
\right.
\end{equation}
and $\sigma_0, \sigma_1, \sigma_2, \sigma_3$ are defined in (\ref{e*}).

Combining all previous arguments we get the following.

\begin{proposition}\label{P1}
The following assertions hold:
\begin{itemize}
\item[1)] 
Relations (\ref{e*}), (\ref{al1}), (\ref{al2}) define $G_0$-grading on the matrix algebra 
$R=M_{2^q}(F)$, where $G_0$ is the elementary abelian 2-group defined in (\ref{al0});
\item[2)]
$\dim R_g=1$ for every $g\in G_0$;
\item[3)]
$R$ has a homogeneous in $G_0$-grading basis consisting of products (\ref{e**}) and any basis 
element is either symmetric or skew-symmetric under transpose involution;
\item[4)]
Every non-zero homogeneous element is invertible; and
\item[5)]
Lie subalgebra $sl_{2^q}$ of traceless matrices is homogeneous in this grading.
\end{itemize}
\end{proposition}
\hfill $\Box$

Applying Proposition \ref{P1}, we construct a grading on some simple Lie superalgebras.
Recall that $P(t)$ (in the notation \cite{Kac}) is a Lie superalgebra $L\subset M_{t,t}(F)$
with
$$
L^{(0)} = \left\{ \left(
           \begin{array}{cc}
             A & 0 \\
             0 & -A^T \\
           \end{array}
         \right)\right\},\,\quad
               L^{(1)} = \left\{ \left(
           \begin{array}{cc}
             0 & B \\
             C & 0 \\
           \end{array}
         \right)
\right\}
$$
where $A,B$ and $C$ are $t\times t$ matrices, $tr A=0$, $B^T=B, C^T=-C$ and $X\to X^T$ is the
transpose involution on $M_t(F)$. We equip $L$ with an abelian grading in the following way.
Let $$t=2^q,  \ \ R=R_1\otimes\cdots\otimes R_q, \ \  R_1=\cdots =R_q= M_2(F)$$ and let $G_0$ be as in
(\ref{al0}). We extend $G_0$ to 
$$G=<a_0>_2\times G_0\simeq (\mathbb Z_2)^{2q+1}$$ 
and define 
$G$-grading on $L$ compatible with canonical $\mathbb Z_2$-grading. If $X_g\in R$ is homogeneous,
$\deg X_g=g\in G_0$, then
\begin{equation}\label{bet1}
Y = \left\{ \left(
           \begin{array}{cc}
             X_g & 0 \\
             0 & -X^T_g \\
           \end{array}
         \right)\right\},
\end{equation}
is homogeneous in $L$, $\deg Y=g$ for all $X_g\in sl_{2^q}(F)\subset R$,

\begin{equation}\label{bet2}
\hbox{if}\quad  X_g \quad\hbox{is symmetric then}\quad
Y = \left\{ \left(
           \begin{array}{cc}
             0 & X_g \\
             0 & 0 \\
           \end{array}
         \right)\right\}
\end{equation}
is homogeneous, $\deg Y= a_0g$

\begin{equation}\label{bet3}
\hbox{if}\quad  X_g \quad\hbox{is skew then}\quad
Y = \left\{ \left(
           \begin{array}{cc}
             0 & 0 \\
             X_g & 0 \\
           \end{array}
         \right)\right\}
\end{equation}
is homogeneous, $\deg Y= a_0g$. 
The following proposition is an immediate consequence of Proposition \ref{P1} and multiplication
rule of $L$.
\begin{proposition}\label{P2}
Let $$G_0=<a_1>_2\times<b_1>_2\times\cdots\times <a_q>_2\times<b_q>_2$$ and $$G=<a_0>_2\times G_0$$ be
elementary abelian 2-groups. Then (\ref{bet1}), (\ref{bet2}) and (\ref{bet3}) define a $G$-grading
on $L=P(2^q)$ compatible with the canonical $\mathbb Z_2$-grading.
 All homogeneous components of $L$
are 1-dimensional. 
If 
$$g=a_0g_0, h=a_0h_0,g_0,h_0\in G_0, \ \  0\ne X_g\in L_g, X_h\in L_h$$
 and both 
$X_g, X_h$ are either of the type (\ref{bet2}) or of the type (\ref{bet3}) then $[X_g,X_h]=0$.
In all other cases $[X_g,X_h]$ is an invertible element of $M_{2^q}(F)$.\hfill $\Box$
\end{proposition}

\section{Graded PI-exponent}

We recall some key notions from the theory of  identities and their numerical invariants. 
We refer the reader to 
 \cite{BahtB, DRbook,GZbook} for details. Consider an absolutely free algebra
$F\{X\}$ with a free generating set
$$
X=\bigcup_{g\in G} X_g,\quad |X_g|=\infty \quad \mbox{for any} \quad g\in G.
$$
One can define a $G$-grading on $F\{X\}$ by setting $\deg_G x=g$, when $x\in X_g$, and
extend this grading to the entire $F\{X\}$ in the natural way. A polynomial $f(x_1,\ldots, x_n)$
in homogeneous variables $x_1\in X_{g_1},\ldots, x_n\in X_{g_n}$ is called a {\it graded identity}
of a $G$-graded algebra $A$ if $f(a_1,\ldots, a_n)=0$ for any  
$a_1\in A_{g_1},\ldots, a_n\in A_{g_n}$. The set $Id^{G}(A)$ of all graded identities of $A$
forms an ideal of $F\{X\}$ which is stable under graded homomorphisms $F\{X\}\rightarrow F\{X\}$.

First, let $G$ be finite, $G=\{g_1,\ldots,g_k\}$ and 
$$X=X_{g_1}\bigcup\ldots\bigcup X_{g_k}.$$ Denote by $P_{n_1,\ldots,n_k}$ the subspace of  $F\{X\}$
of multilinear polynomials of total degree $n=n_1+\cdots+n_k$ in variables
$$x^{(1)}_{1},\ldots,x^{(1)}_{n_1}\in X_{g_1},\ldots, x^{(k)}_{1},\ldots,x^{(k)}_{n_k}\in X_{g_k}.$$
Then the value
$$
c_{n_1,\ldots,n_k}(A)=\dim\frac{P_{n_1,\ldots,n_k}}{P_{n_1,\ldots,n_k}\cap Id^{G}(A)}
$$
is called a {\it partial codimension} of $A$ while
\begin{equation}\label{ee1}
c_n^{G}(A)=\sum_{n_1+\cdots+n_k=n} {n\choose n_1,\ldots,n_k} c_{n_1,\ldots,n_k}(A)
\end{equation}
is called a {\it graded codimension} of $A$. Recall that the {\it support of the grading} is the set
$$
Supp~A=\{g\in G\vert A_g\ne 0\}.
$$
Note that if $Supp~A\ne G$, say,  $Supp~A=\{g_1,\ldots,g_d\}$, $d<k$, then the value
\begin{equation}%\label{ee2}
\sum_{n_1+\cdots+n_d=n} {n\choose n_1,\ldots,n_d} \dim\frac{P_{n_1,\ldots,n_d}}{P_{n_1,\ldots,n_d}\cap Id^{G}(A)}
\end{equation}
coincides with (\ref{ee1}). 

Denote
\begin{equation}\label{ee3}
P_{n_1,\ldots,n_k}(A)=\frac{P_{n_1,\ldots,n_k}}{P_{n_1,\ldots,n_k}\cap Id^{G}(A)}.
\end{equation}

For finding a lower bound for PI-exponent we need the following observation.

\begin{lemma}\label{L1}
Let $A$ be a $G$-graded algebra with the support $Supp A=\{g_1,\ldots,g_d\}\subseteq G$. Let also 
$\dim A_g=1$ for any $g\in Supp A$. Then
\begin{itemize}
\item[(1)]
if $P_{n_1,\ldots,n_d}(A)\ne 0$ then $\dim P_{n_1,\ldots,n_d}(A)=1$,
\item[(2)]
$\dim P_{n_1,\ldots,n_d}(A)=1$ if and only if there exist $u_1\in A_{g_1},\ldots,u_d\in A_{g_d}$ and a
monomial $m(u_1,\ldots,u_d)=m\ne 0$ on $u_1,\ldots,u_d$ such that every $u_j$ 
appears in $m$ exactly $n_j$
times, $j=1,\ldots, d$.
\end{itemize}
\end{lemma}

{\it Proof}. First, let $P_{n_1,\ldots, n_d}(A)\ne 0$. Then there exists a multilinear homogeneous polynomial
$$f=f(x_1^{(1)},\ldots, x_{n_1}^{(1)},\ldots, x_{1}^{(d)},\ldots, x_{n_d}^{(d)})\in P_{n_1,\ldots,n_d}$$
which is not an identity of $A$. That is, one can find $u_1\in A_{g_1},\ldots,u_d\in A_{g_d}$ such that
$f(u_1,\ldots, u_d)\ne 0$. If
$$g=g(x_1^{(1)},\ldots, x_{n_1}^{(1)},\ldots, x_{1}^{(d)},\ldots, x_{n_d}^{(d)})
\in P_{n_1,\ldots,n_d}\setminus Id^G(A)$$ then
$$
g(u_1,\ldots,u_1,\ldots,u_d,\ldots,u_d)=\lambda f(u_1,\ldots,u_1,\ldots,u_d,\ldots,u_d)
$$
for some scalar $\lambda$ since $\dim A_g=1$ for $g=g_1^{n_1}\cdots g_d^{n_d}$. Hence $g-\lambda f\equiv 0$
is an identity of $A$. This proves (1).

Now let $\dim P_{n_1,\ldots,n_d}(A)=1$, that is $P_{n_1,\ldots,n_d}(A)\ne 0$. Then there exist
$$f=f(x_1^{(1)},\ldots, x_{n_1}^{(1)},\ldots, x_{1}^{(d)},\ldots, x_{n_d}^{(d)})\in P_{n_1,\ldots,n_d}\setminus
Id^G(A)$$ and $u_1\in A_{g_1},\ldots,u_d\in A_{g_d}$ such that
$$
f(\underbrace{u_1,\ldots,u_1}_{n_1},\ldots,\underbrace{u_d,\ldots,u_d}_{n_d})\ne 0
$$
in $A$. Hence, at least one monomial of $f$ has a non-zero value under evaluation
$\varphi: F\{X\}\mapsto A,$ where 
$$\varphi(x_j^{(i)})=u_i, \ \ 1\le i\le d, \ \  1\le j\le n_i.$$ This implies $(2)$,
and have we completed the proof.

\hfill $\Box$

\begin{corollary}\label{C1}
\begin{equation}\label{eqq1}
c_n^G=\sum {n\choose n_1,\ldots, n_d}
\end{equation}
where the sum in (\ref{eqq1}) is taken over all tuples $(n_1,\ldots, n_d)$ such that
\begin{equation}\label{eqq2}
P_{n_1,\ldots, n_d}(A)\ne 0.
\end{equation}
Moreover, for the inequality (\ref{eqq2}) it suffices to check the condition (2) of Lemma \ref{L1}.
\end{corollary}
\hfill $\Box$

Now we go back to the Lie superalgebra $$L=L^{(0)}\oplus L^{(1)}=P(t),  \ \ t=2^q,$$ with the
$G$-grading presented in Proposition \ref{P2}. First, we give an upper bound for $exp^G(L)$. 
Note that Stirling formula for factorials implies the inequalities
\begin{equation}\label{eq31}
\frac{1}{n^d}\Phi(n;n_1,\ldots,n_d)^n \le {n\choose n_1,\ldots,n_d}\le n \Phi(n;n_1,\ldots,n_d)^n 
\end{equation}
where
$$
\Phi(n;n_1,\ldots,n_d)=(\frac{n_1}{n})^{-\frac{n_1}{n}}\cdots (\frac{n_d}{n})^{-\frac{n_d}{n}}
$$
and $n=n_1+\cdots+n_d$.

Denote
$$
a=\frac{t(t+1)}{2}, b=\frac{t(t-1)}{2}, c=t^2-1, d=a+b+c=\dim L.
$$
The algebra $L$ has a natural $\mathbb Z$-grading
$$
L=\mathcal L_{-1}\oplus \mathcal L_{0}\oplus \mathcal L_{1}
$$
where
$$
\mathcal L_{-1}=
 \left\{ \left(
           \begin{array}{cc}
             0 & 0 \\
             C & 0 \\
           \end{array}
         \right)  \right\},\,
\mathcal L_{0}=L^{(0)}=
 \left\{ \left(
           \begin{array}{cc}
             A & 0 \\
             0 & -A^T \\
           \end{array}
         \right)  \right\},\,
\mathcal L_{1}=
 \left\{ \left(
           \begin{array}{cc}
             0 & B \\
             0 & 0 \\
           \end{array}
         \right)  \right\}.
$$
All remaining components $\mathcal L_{k}$, $k\ne0,\pm 1$, are zero. Clearly,
$P_{n_1,\ldots,n_d}(L)\ne 0$ only if
\begin{equation}\label{eq4}
|n_1+\cdots+n_a-n_{a+1}\cdots-n_{a+b}|\le 1
\end{equation}
where $\{g_1,\ldots,g_d\}\subseteq G$ is the support $Supp L$. It follows from Corollary \ref{C1} 
and (\ref{eq31})  that
\begin{equation}\label{eq5}
\frac{1}{n^d}
\hbox{\rm max}
\{\Phi(n;n_1,\ldots,n_d)^n\}
 \le c_n^G(L)\le n^d 
 \hbox{\rm max}
\{ \Phi(n;n_1,\ldots,n_d)^n\}
\end{equation}
where the maximum is taken among all $n_1,\ldots,n_d$ satisfying (\ref{eq4}).

First, consider the case where the left  side of (\ref{eq4}) is equal to zero. Then we rewrite
$$
\Phi(n;n_1,\ldots,n_d)=\Phi(x_1,\ldots,x_d)
$$
where $x_1+\cdots+x_d=1$, $x_1,\ldots,x_d\ge 0$,
\begin{equation}\label{eq6}
\Phi(x_1,\ldots,x_d)=x_1^{-x_1}\cdots x_d^{-x_d}
\end{equation}
and $$x_1+\cdots +x_a= x_{a+1}+\cdots+x_{a+b}.$$ It is easy to see that the
maximal value of the function (\ref{eq6})  is achieved when
$$
x_1=\cdots =x_a,\, x_{a+1}=\cdots=x_{a+b},\, x_{a+b+1}=\cdots=x_{a+b+c}.
$$
Denote $x=x_1, y=x_{a+b}, z=x_{a+b+c}$. Then (\ref{eq6}) does not exceed
$$
\widetilde\Phi=\widetilde\Phi(x,y,z)=x^{-ax}y^{-by}z^{-cz}
$$ 
and $x,y,z$ satisfy the relations $ax=by$, $ax+by+cz=1$. These relations
imply
$$
\widetilde\Phi^{-1} = z^{(t^2-1)z}(1-(t^2-1)z)^{(1-(t^2-1)z)}
(t^2(t^2-1))^\frac{(t^2-1)z-1}{2}
$$
as a function of $z$. Then
$$
g(z)=\ln \widetilde\Phi^{-1}= cz\ln z+(1-cz)\ln(1-cz)-\frac{1}{2}(1-cz)\ln(ct^2).
$$
Direct calculations show that $g'(z)=0$ only if
$$
z=z_0=(t^2-1+t\sqrt{t^2-1})^{-1}
$$
and $g''(z_0)>0$. Hence, in $z_0$ the function $g(z)$ has a local mnimum. Moreover,
$$
g(z_0)=-\ln (t^2-1+t\sqrt{t^2-1}).
$$
It follows that
$$
\widetilde\Phi\le t^2-1+t\sqrt{t^2-1}
$$
and
\begin{equation}\label{eq7}
\sqrt[n]{c_n^G(L)}\le n^\frac{d}{n}(t^2-1+t\sqrt{t^2-1}) 
\end{equation}
as follows from (\ref{eq5}) in the case $n_1+\cdots+n_a=n_{a+1}+\cdots+n_{a+b}$.

If $n_1+\cdots+n_a-n_{a+1}-\cdots-n_{a+b}=-1$ then
$$
{n\choose n_1,\ldots,n_d} \le {n+1\choose n_1+1,n_2,\ldots,n_d}
$$
and
\begin{equation}\label{eq8}
\sqrt[n]{c_n^G(L)}\le (n+1)^\frac{d}{n+1}(t^2-1+t\sqrt{t^2-1}).
\end{equation}
Similarly, if  $n_1+\cdots+n_a-n_{a+1}-\cdots-n_{a+b}=1$ then
\begin{equation}\label{eq9}
\sqrt[n]{c_n^G(L)}\le (n-1)^\frac{d}{n-1}(t^2-1+t\sqrt{t^2-1})
\end{equation}
since
$$
{n\choose n_1,\ldots,n_d} \le n{n-1\choose n_1,\ldots,n_{a+b-1}, n_{a+b}-1,n_{a+b+1},\ldots,n_d}.
$$

Inequalities (\ref{eq7}), (\ref{eq8}) and (\ref{eq9}) give us the following.
\begin{lemma}\label{L2}
$$ 
exp^G(L)\le t^2-1 +t\sqrt{t^2-1}.
$$
\end{lemma}

Now we will get the same lower bound.

\begin{lemma}\label{L3}
\begin{equation}\label{lb}
exp^G(L)\ge t^2-1 +t\sqrt{t^2-1}.
\end{equation}
\end{lemma}

{\it Proof}. Recall that $L$ is $\mathbb Z$-graded algebra, 
$L=\mathcal L_{-1}\oplus \mathcal L_{0}\oplus \mathcal L_{1}$, and $a=\dim\mathcal L_1$,
$b=\dim\mathcal L_{-1}$, $c=\dim\mathcal L_0$. Consider a collection
$$
X=\{\underbrace{x_1,\ldots,x_1}_{b},\ldots, \underbrace{x_a,\ldots,x_a}_{b}\}
$$
where $x_1,\ldots, x_a$ are homogeneous in $G$-grading elements $\mathcal L_1$ with
pairwise distinct degree in $G$-grading. Similarly, we take
$$
Y=\{\underbrace{y_1,\ldots,y_1}_{a},\ldots, \underbrace{y_b,\ldots,y_b}_{a}\},
$$
with homogeneous $y_1,\ldots, y_b\in \mathcal L_{-1}$, $deg_G y_i$ are distinct. Renaming
elements of $X,Y$ we write
$$
X=\{x^{(1)},\ldots,x^{(ab)}\},\quad Y=\{y^{(1)},\ldots,y^{(ab)}\}.
$$
We remark that any $x_i, 1\le i\le a$  appears among $x^{(1)},\ldots,x^{(ab)}$ exactly
$b$ times. Similarly, any $y_j, 1\le j\le  b$, appears among $y^{(1)},\ldots,y^{(ab)}$
exactly $a$ times. Consider supercommutators
$$
z_1=[x^{(1)}, y^{(1)}],\ldots, z_{ab}=[x^{(ab)}, y^{(ab)}].
$$
By Proposition \ref{P2} all $z_i$ are invertible in $M_{2t}(F)$ matrices homogeneous in 
$G$-grading of $L$. Also, $$z_1,\ldots,z_{ab}\in L^{(0)}\simeq sl_{2t}(F).$$ Note that $xy=\pm yx$ 
for any homogeneous $x,y\in L^{(0)}$. It follows that for any $i=1,\ldots, ab$ there
exists $z_i'\in L^{(0)}$ homogeneous in $G$-grading such that
$$
[z_i',z_i]=2z_i' z_i\ne 0
$$
where the product $z_i' z_i$ is taken in the associative algebra $M_{2t}(F)$.
Hence, the left-normed Lie commutators
$$
z_k^{(i)}= [z_i', \underbrace{z_i,\ldots,z_i}_{k}]=2^kz_i' z_i^k, \, k=1,2,\ldots,
$$
are non-zero homogeneous elements of $L^{(0)}$.

As before, one can find homogeneous $u_1,\ldots, u_{ab}\in L^{(0)}$ and linearly independent
homogeneous $v_1,\ldots, v_{c}\in L^{(0)}$ such that
$$
w_k=[z_k^{(1)},u_1,z_k^{(2)},u_2,\ldots, z_k^{(ab)},u_{ab}]\ne 0
$$
and 
$$
w_{k,s}=[w_k,w_1',\underbrace{v_1,\ldots,v_1}_{s},w_2',\underbrace{v_2,\ldots,v_2}_{s},
\ldots, w_c',\underbrace{v_c,\ldots,v_c}_{s}]\ne 0
$$
for some homogeneous $w_1',\ldots,w_c'\in L^{(0)}$.

If $u$ is a monomial on $x_1,\ldots,x_a,y_1,\ldots,y_b,v_1,\ldots,v_c$ in $L$ then we will denote by
${\rm Deg}_{x_i}u, {\rm Deg}_{y_i}u, {\rm Deg}_{v_i}u$ the total number of factors $x_i,y_i$
and $v_i$ in $u$, respectively. Then
$$
{\rm Deg}_{x_i}w_{k,s} \ge kb \quad\hbox{for all}\quad i=1,\ldots,a,
$$
$$
{\rm Deg}_{y_i}w_{k,s} \ge ka \quad\hbox{for all}\quad i=1,\ldots,b,
$$
$$
{\rm Deg}_{v_i}w_{k,s} \ge s \quad\hbox{for all}\quad i=1,\ldots,c.
$$
Total degrees ${\rm Deg}$ on $\{x_\alpha,y_\beta,v_\gamma \}$ are as follows:
$$
{\rm Deg} z_k^{(i)}=2k+1, {\rm Deg} w_k=2ab(k+1), {\rm Deg} w_{k,s}=2ab(k+1)+c(s+1)=n.
$$
Denote
\begin{equation}\label{pp1}
n_i={\rm Deg}_{x_i} w_{k,s},\, i=1,\ldots, a,
\end{equation}
\begin{equation}\label{pp2}
n_{a+i}={\rm Deg}_{y_i} w_{k,s},\, i=1,\ldots, b,
\end{equation}
\begin{equation}\label{p3}
n_{a+b+i}={\rm Deg}_{v_i} w_{k,s},\, i=1,\ldots, c.
\end{equation}
If $$m_1=\cdots=m_a=kb, m_{a+1}=\cdots=m_{a+b}=ka, m_{a+b+1}=\cdots=m_{a+b+c}=s,$$
and
$$m=m_1+\cdots+m_{a+b+c}=2abk+cs$$ then $n-m=2ab+c$ and
\begin{equation}\label{eq10}
{n\choose n_1,\ldots, n_d} \ge {m\choose m_1,\ldots, m_d} \ge$$ $$\ge\frac{1}{m^d}
\Phi(m;m_1,\ldots, m_d)^m= \frac{1}{m^d}\widetilde \Phi(\frac{kb}{m},\frac{ka}{m},\frac{s}{m})^m.
\end{equation}

Denote $\frac{k}{s}=\alpha$. Then
$$
\frac{s}{m}=\frac{s}{2abk+cs}=\frac{1}{\frac{k}{s}\cdot\frac{t^2(t^2-1)}{2}+t^2-1}=
\frac{1}{t^2-1+\alpha\frac{t^2(t^2-1)}{2}}.
$$

Note that if $$\beta=\frac{2}{t\sqrt{t^2-1}}$$ then
$$
t^2-1+\beta\frac{t^2(t^2-1)}{2})=t^2-1+t\sqrt{t^2-1}
$$
and
$$
\widetilde\Phi(\bar x,\bar y,\bar z)=\Phi_{max}=t^2-1+t\sqrt{t^2-1}
$$
provided that 
$$
\bar z = (t^2-1+t\sqrt{t^2-1})^{-1},\, a\bar x=b\bar y,\, a\bar x+d\bar y+c\bar z=1.
$$
In particular, if $$\alpha=\frac{k}{s}\to\beta$$ then
$$
\widetilde \Phi(\frac{kb}{m},\frac{ka}{m},\frac{s}{m})\to\Phi_{max}.
$$
More precisely, for any $\varepsilon>0$ there exists real $\delta$ such that the inequality
\begin{equation}\label{eq11}
|\frac{k}{s}-\frac{2}{t\sqrt{t^2-1}}| < \delta
\end{equation}
implies
\begin{equation}\label{eq12}
\Phi(m;m_1,\ldots, m_d) \ge t^2-1+t\sqrt{t^2-1}-\varepsilon.
\end{equation}

Fix one pair $(k,s)$ with the relation (\ref{eq11}) and take 
$$m=2abk+cs, \ \  \bar n_1=m+2ab+c, \ \ 
n_i$$ as in (\ref{pp1}), (\ref{pp2}), (\ref{p3}). Then we have for any $r=1,2,\ldots$,
$$
{r\bar n_1\choose rn_1,\ldots, rn_d} \ge \frac{1}{(r\bar n_1)^d} \Phi(rm;rm_1,\ldots,rm_d)^{rm}
=\frac{1}{(r\bar n_1)^d} \Phi(m;m_1,\ldots,m_d)^{rm}
$$
$$ \ge
\frac{1}{(r\bar n_1)^d}(t^2-1+t\sqrt{t^2-1}-\varepsilon)^{rm}
$$
as follows from (\ref{eq10}), (\ref{eq12}).

Denote $\bar n_r=r\bar n_1$. For any given $\rho>0$ we can choose $\bar n_1$ large enough
and suppose that
$$
\frac{rm}{\bar n_r}=\frac{\bar n_r-(2ab+c)r}{\bar n_r}=1-\frac{2ab+c}{\bar n_1} > 1-\rho
$$
from which it follows that
\begin{equation}\label{eq13}
c_{\bar n_r}^G \ge\frac{1}{\bar n_r^d}(t^2-1+t\sqrt{t^2-1}-\varepsilon)^{1-\rho}.
\end{equation}
Since $\bar n_{r+1}-\bar n_r=2ab+c=const$ and
$$
{n'\choose n_1',\ldots,n_d'} \ge {n\choose n_1,\ldots,n_d}
$$
as soon as 
$$n'=n+1, n_1'\ge n_1,\ldots, n_d'\ge n_d,$$
  (\ref{eq13}) implies the inequality
$$
exp^G(L)\ge t^2-1+t\sqrt{t^2-1}-\varepsilon.
$$
Recall that $\varepsilon >0$ is arbitrary, hence (\ref{lb}) follows and we are done.
\hfill $\Box$

\bigskip

{\it Proof of Theorem \ref{T1}.} The assertions of Theorem \ref{T1} now follow  from Lemmas \ref{L2} and \ref{L3}.
\hfill $\Box$

\end{document}